\documentclass[11pt]{amsart}
\usepackage{amsmath}
\usepackage{graphicx}
\usepackage{amssymb}

\newtheorem{thm}{Theorem}[section]
\newtheorem{lem}[thm]{Lemma}

\newtheorem{cor}[thm]{Corollary}

\theoremstyle{definition}

\newtheorem{exmp}[thm]{Example}

\renewcommand{\bold}[1]{\emph{#1}}
\newcommand{\orb}[1]{\mathcal{O}_{#1}}
\newcommand{\Fset}[1]{\mathcal{F_{#1}}}
\newcommand{\Wset}[1]{\mathcal{W_{#1}}}
\newcommand{\Dset}[1]{\mathcal{D_{#1}}}

\newcommand{\trans}[2]{\overset{\t_{#1}}{\rule{0mm}{3.5mm}(1\,#2)}}

\newcommand{\s}{{\sigma}}

\newcommand{\g}{{\gamma}}
\newcommand{\p}{{\pi}}

\renewcommand{\t}{{\tau}}

\newcommand{\map}[3]{{#1} \! : \! {#2} \! \longrightarrow \! {#3}}

\newcommand{\wo}{\setminus}

\newcommand{\Sym}[1]{{{\mathfrak{S}_{#1}}}}

\begin{document}



\title{Minimal Factorizations of Permutations Into Star Transpositions}


\author{J. Irving}
\address{Department of Mathematics and Computing Science \\
        St. Mary's University, Halifax, NS, B3H 3C3, Canada}
\author{A. Rattan}
\address{Department of Mathematics \\
        Massachusetts Institute of Technology,
    Cambridge, MA, 02139, USA}

\begin{abstract}
We give a compact expression for the number of factorizations of any permutation into a minimal number of
transpositions of the form $(1\,i)$.  This generalizes earlier work of Pak in which substantial restrictions were
placed on the permutation being factored. Our result exhibits an unexpected and simple symmetry of star
factorizations that has yet to be explained in a satisfactory manner.
\end{abstract}

\maketitle

\section{Introduction}

It is well known that the symmetric group $\Sym{n}$ is generated by various sets of
transpositions, and it is natural  to ask for the number of decompositions of a
permutation into a minimal number of factors from such a set.  For instance, a famous
paper of D\'enes~\cite{denes} addresses this question when the generating set is
taken to consist of all transpositions.  Stanley~\cite{stanley-reduced} has also
considered the problem for the set of Coxeter generators $\{(i\,\,i+1) \,:\, 1 \leq i
< n\}$.

More recently, Pak~\cite{pak} considered minimal decompositions of permutations
relative to the generating set $S=\{(1\,i) \,:\, 2 \leq i \leq n\}$.  The elements of
$S$ are called \bold{star transpositions} because the labelled graph on vertex set
$[n]=\{1,\ldots,n\}$ obtained from them by interpreting $(a\,b)$ as an edge between
vertices $a$ and $b$ is star-shaped. Pak proves that any permutation $\p \in \Sym{n}$
that fixes 1 and has $m$ cycles of length $k \geq 2$ admits exactly
\begin{equation}
\label{eq:pak}
        \frac{k^{m} (mk+m)!}{n!}
\end{equation}
decompositions into the minimal number $n+m-1$ of star transpositions. He leaves open the problem of
extending~\eqref{eq:pak} to more general target permutations $\p$, and it is the purpose of this paper to answer
this question.

Our result is best expressed in terms of \emph{minimal transitive star factorizations}, which we now define. A
\bold{star factorization} of $\p \in \Sym{n}$ of \emph{length} $r$ is an ordered list $f=(\t_1,\ldots,\t_r)$ of
star transpositions $\t_i$ such that $\t_1 \cdots \t_r = \p$.\footnote{We multiply permutations in the usual
order, so $\rho\s(j)=\rho(\s(j))$.}  We say $f$ is \bold{minimal} if $\p$ admits no star factorization of length
less than $r$, and \bold{transitive} if the group generated by its factors acts transitively on $[n]$.

Observe that a permutation $\pi=(1\,\, b_2 \cdots b_{\ell_1}) (a^2_1\cdots a^2_{\ell_2}) \cdots (a^m_1 \cdots
a^m_{\ell_m}) \in \Sym{n}$ with $m$ cycles admits the transitive star factorization
\begin{equation*}
    \begin{split}
        \pi = \underbrace{(1\,\, b_{\ell_1})
        (1\,\, b_{\ell_1-1})\cdots
        (1\,\, b_2)}_{= (1\,\, b_2 \cdots b_{\ell_1})}
        \underbrace{(1\,\, a^2_1)
        (1\,\, a^2_{\ell_2})
        (1\,\, a^2_{\ell_2 -1})\cdots
        (1\,\, a^2_1)}_{= (a^2_1 \cdots a^2_{\ell_2})}\cdots\\
        \underbrace{(1\,\, a^m_1)
        (1\,\, a^m_{\ell_m})
        (1\,\, a^m_{\ell_m -1})\cdots
        (1\,\, a^m_1)}_{= (a^m_1 \cdots a^m_{\ell_m})}
    \end{split}
\end{equation*}
of length $\ell_1 - 1 + \sum_{i=2}^m (\ell_i + 1) = n+m-2$. Moreover, it is well known~\cite[Proposition
2.1]{gj-transitive} that \emph{any} transitive star factorization of $\p$ requires at least this many
factors.\footnote{In fact, this holds true when arbitrary transposition factors are allowed.} Thus a transitive
star factorization of $\p$ of length exactly $n+m-2$ is said to be \emph{minimal transitive}.

Our main result is the following:
\begin{thm}
\label{thm:mainthm} Let $\p \in \Sym{n}$ be any permutation with cycles of lengths $\ell_1,\ldots,\ell_m$. Then
there are precisely
$$
        \frac{(n+m-2)!}{n!} \ell_1 \cdots \ell_m
$$
minimal transitive star factorizations of $\p$.
\end{thm}
Notice that Pak's formula~\eqref{eq:pak} is recovered from Theorem~\ref{thm:mainthm}
by setting $\ell_1=1$ and $\ell_2=\cdots=\ell_{m+1}=k$ and observing that a star
factorization of a permutation with no fixed points other than (possibly) 1 must be
transitive, since $\p(a) \neq a$ means any star factorization of $\p$  involves the
factor $(1\,a)$.

Given the special role played by the symbol 1 in star factorizations, the lack of bias towards this symbol in
the enumerative formula of Theorem~\ref{thm:mainthm} is quite surprising.  Indeed, this symmetry is a very
compelling aspect of the theorem, and it is not yet understood.

Consider now a permutation $\p$ having fixed points $i_1,\ldots,i_k$ and possibly 1.
A minimal (though not transitive) star factorization of $\p$ should certainly
\emph{not} contain any of the factors $(1\,i_1),\ldots,(1\,i_k)$.  Indeed, $\p$
naturally induces a permutation $\p'$ on $[n] \wo \{i_1,\ldots,i_k\}$ having no fixed
points other than (possibly) 1, and minimal star factorizations of $\p$ are simply
minimal transitive star factorizations of $\p'$.  Since $\p'$ has $m-k$ cycles when
$\p$ has $m$ cycles, we obtain the following result by setting $n=n-k$ and $m=m-k$ in
Theorem~\ref{thm:mainthm}.
\begin{cor}
\label{thm:maincor} Let $\p \in \Sym{n}$ be any permutation with cycles of lengths
$\ell_1,\ldots,\ell_m$ including exactly $k$ fixed points not equal to 1. Then there
are
$$
        \frac{(n+m-2(k+1))!}{(n-k)!} \ell_1 \cdots \ell_m
$$
minimal  star factorizations of $\p$.
\end{cor}

We prove Theorem~\ref{thm:mainthm} in two stages. In Section~\ref{sec:characterization}, we begin by giving a
complete characterization of minimal transitive star factorizations (Lemma~\ref{lem:characterize}). We then use
this characterization in Section~\ref{sec:counting} to build a correspondence between star factorizations and
certain restricted words, finally using the cycle lemma to count these words and hence prove
Theorem~\ref{thm:mainthm}.

This path to Theorem~\ref{thm:mainthm} is deliberately similar to that followed in~\cite{pak}. However, in
Section~\ref{sec:graphical}, we briefly describe an elegant graphical approach to this problem that employs
the well-known connection between factorizations of permutations and embeddings of graphs on surfaces
(\emph{i.e.} maps).  Finally, Section~\ref{sec:final} contains some brief comments on recent extensions of
Theorem~\ref{thm:mainthm} and its curious symmetry.

\section{Characterizing Star Factorizations}
\label{sec:characterization}

\emph{Throughout this section we have in mind a fixed permutation $\p \in \Sym{n}$ and a minimal transitive
star factorization $f=(\t_1,\ldots,\t_r)$ of $\p$.}

Our arguments are  best understood with a concrete example at hand. For this purpose, we will often refer to
the factorization
\begin{equation}
\label{eq:fexmp}
        \trans{1}{9}
        \trans{2}{\,11}
        \trans{3}{9}
        \trans{4}{2}
        \trans{5}{\,10}
        \trans{6}{5}
        \trans{7}{3}
        \trans{8}{3}
        \trans{9}{4}
        \trans{10}{7}
        \trans{11}{6}
        \trans{12}{6}
        \trans{13}{\,10}
        \trans{14}{8}
\end{equation}
of
\begin{equation}
\label{eq:pexmp}
    \p =    (1\,\,8\,\,2)
            (3)
            (4\,\,5\,\,10\,\,7)
            (6)
            (9\,\,11) \in \Sym{11}.
\end{equation}

Let us say that a transposition $(1\,i)$ \bold{meets} a cycle $\s$ (and vice versa) if $\s$ contains the
symbol $i$. We say a factor $\t_i$ is to the \emph{left}  of $\t_j$ if $i < j$, and to the \emph{right} if $j
> i$. So, in the example above, $(1\,5)$ meets $(4\,\,5\,\,10\,\,7)$, and $(1\,2)$ is to the left of
$(1\,7)$. Clearly every star transposition meets exactly one cycle of $\p$.

Our goal here is to characterize minimal transitive star factorizations of $\p$.

\begin{lem}
\label{lem:order} Let $\s$ be a cycle of $\p$.
\begin{enumerate}
\item   If $\s = (a_1\,a_2\,\ldots\,a_{\ell})$, where $a_i \neq 1$ for all $i$,
    then some transposition $(1\,a_j)$ appears exactly twice in $f$, while all
    transpositions $(1\,a_i)$ with $i \neq j$ appear exactly once.

\item   If $\s = (1\,b_2\,\cdots\,b_{\ell})$,  then
        each transposition $(1\,b_i)$ appears only once in $f$.
\end{enumerate}
Moreover, in the first case, if  $(1\,a_1)$ appears twice, then the factors of $f$
meeting $\s$ appear in left-to-right order
$(1\,a_{1}), (1\,a_\ell),\ldots, (1\,a_{2}), (1\,a_1)$.  In the second case, the
factors meeting $\s$ appear in left-to-right order
$(1\,b_{\ell}), (1\,b_{\ell-1}), \ldots, (1\,b_{2})$.
\end{lem}

\begin{proof}
Suppose $\s = (a_1\,a_2\,\ldots\,a_{\ell})$ with  $a_i \neq 1$.  It is clear that for
$f$ to be transitive every transposition $(1\,a_i)$ must appear at least once as a
factor.  Let $(1\,a_j)$ be the leftmost (last in order of multiplication) factor of
$f$ that meets $\s$. If $(1\,a_j)$  appeared only this once, then we would have $\p =
\p_1\, (1\,a_j)\,\p_0$, where $\p_0$ fixes $a_j$ and $\p_1$ fixes all $a_i$. In
particular, $\s(a_j)=\p(a_j) = \p_1(1) \neq a_i$ for any $i$, a contradiction.

On the other hand, if $\s = (1\,b_2\,\cdots\,b_{\ell})$ then  again transitivity requires that $f$ contain factors
$(1\,b_2), \ldots, (1\,b_{\ell})$.  So if the cycles of $\p$ are $\s_1,\ldots,\s_m$, where $\s_1$ contains symbol
1, then $\s_i$ meets at least $\ell_i+1$ factors of $f$ for $i \neq 1$, while $\s_1$ meets at least $\ell_1-1$
factors. That is, $f$ has at least $(\ell_1-1) + \sum_{i=2}^m (\ell_i+1) = n+m-2$ factors. But since $f$ is
minimal transitive, it has exactly this many factors. Hence all factors are accounted for and parts (1) and (2) of
the lemma follow. It remains to determine the relative ordering of the factors meeting $\s$.

We return to the case $\s = (a_1\,a_2\,\ldots\,a_{\ell})$ with  $a_i \neq 1$, and assume without loss of
generality that $(1\,a_1)$ appears twice in $f$.   The proof given above identified $(1\,a_1)$ as the leftmost
factor of $f$ meeting $\s$. However, reading the factors of $f$ in reverse order yields a factorization $f'$ of
$\p^{-1}$, and the same logic now identifies $(1\,a_1)$ as the leftmost factor of $f'$ meeting $\s$.  Thus
$(1\,a_1)$ appears in $f$ in the leftmost and rightmost positions amongst all factors meeting $\s$.  Finally, note
that for $1 \leq i < \ell$ the factor $(1\,a_{i+1})$ is to the left of the rightmost occurrence of $(1\,a_{i})$ in
$f$, as otherwise we would have $\p = \p_1\,(1\,a_i)\,\p_0$, where $\p_0$ fixes $a_{i}$ and $\p_1$ fixes
$a_{i+1}$, and this gives the contradiction $\s(a_i) = \p(a_i)=\p_1(1) \neq a_{i+1}$.  It follows that the factors
meeting $\s$ appear in order $(1\,a_{1}), (1\,a_{\ell}),\ldots, (1\,a_{2}), (1\,a_1)$.

If instead $\s = (1\,b_2\,\cdots\,b_{\ell})$, then the same logic just applied shows that for $2 \leq i < \ell$,
the factor $(1\,b_{i+1})$ appears to the left of $(1\,b_i)$ in $f$.  Thus the factors meeting $\s$ appear in order
$(1\,b_{\ell}), (1\,b_{\ell-1}), \ldots, (1\,b_{2})$, as claimed. \qed
\end{proof}

The next lemma asserts that the factors of a minimal transitive star factorization are nested in a well defined
manner. This ``non-crossing'' property makes it unsurprising that such factorizations can be encoded as trees.
(See Section~\ref{sec:graphical} for details.)  Note that one immediate consequence of the lemma is that there
exists some cycle of $\p$ such that all the factors meeting this cycle appear consecutively in $f$.

\begin{lem}
\label{lem:between} Let $\s$ and $\hat{\s}$ be distinct cycles of $\p$.  Suppose there exist $s < v < t$ such that
factors $\t_s$ and $\t_t$ of $f$ meet $\s$ while $\t_v$ meets $\hat{\s}$. Then $\hat{\s}$ does not contain the
symbol 1, and all $\t_j$ that meet $\hat{\s}$ have $s < j < t$.
\end{lem}

\begin{proof}
Without loss of generality we may assume that indices $s$ and $t$ are ``extremal'', in the sense that if $\t_{s'}$
and $\t_{t'}$ meet the same cycle of $\p$, with $s < s' < v < t' < t$, then this common cycle is $\hat{\s}$. (If
not, simply restart by letting $\s$ be the common cycle and replacing $s$ and $t$ with $s'$ and $t'$.)

Let $\t_s=(1\,b)$ and $\t_t=(1\,a)$. Since $\t_s$ and $\t_t$ are assumed to meet the same cycle of $\p$,
Lemma~\ref{lem:order} implies $\t_s$ is the leftmost copy of $(1\,b)$ in $f$, $\t_t$ is the rightmost copy of
$(1\,a)$ in $f$, and $\p(a)=b$. It follows from these criteria that the permutation $\t_{s+1} \cdots
\t_{t-1}$ fixes $1$, and therefore $$i := \max\{k \,:\, \text{$k \leq v$ and $\t_{k} \cdots \t_{t-1}$ fixes
1}\}$$ is well defined.  Say $\t_i = (1\,c)$. Notice that this factor must occur twice amongst those of
$\g=\t_i \cdots \t_{t-1}$, as otherwise $\g(c)=1$ and hence $\g$ does not fix 1, contrary to the definition
of $i$.

Suppose $\t_j=\t_i=(1\,c)$ for some $j > i$.  Then Lemma~\ref{lem:order} implies $c$ cannot appear in any factor
between $\t_i$ and $\t_j$, so the permutation $\t_i \cdots \t_j$ fixes 1. But since $\t_i \cdots \t_{t-1}$ fixes
1, it follows that $\t_{j+1} \cdots \t_{t-1}$ also fixes 1.  Thus the maximality of $i$ forces $j \geq v$.
However, if $j > v$ then we have two identical factors $\t_i$ and $\t_j$, with $s < i < v < j < t$, that meet the
same cycle of $\p$, and by hypothesis this common cycle must be $\hat{\s}$.   In this case, Lemma~\ref{lem:order}
rules out the possibility of $\hat{\s}$ containing symbol 1 (because no transposition meeting the cycle containing
1 can appear twice in $f$), and also implies any other factor of $f$ that meets $\hat{\s}$ lies between $\t_i$ and
$\t_j$, as desired. The remaining case is $j=v$, in which $\t_v$ occurs twice between $\t_s$ and $\t_t$. Again the
result follows from Lemma~\ref{lem:order}. \qed
\end{proof}

The statements of Lemmas~\ref{lem:order} and~\ref{lem:between} are crafted with the
implicit assumption that $f=(\t_1,\ldots,\t_r)$ is a minimal transitive star
factorization of $\pi$.  We now show that this can, in fact, be deduced from the
conditions on $f$ established by the lemmas.  That is to say, if $f$ is a star
factorization whose factors are related to the permutation $\p$ in the manner
described by the lemmas above, then $f$ is necessarily a minimal transitive star
factorization of $\p$.

\begin{lem}
\label{lem:characterize} The conditions on $f$ guaranteed by Lemmas~\ref{lem:order}
and~\ref{lem:between} characterize minimal transitive star factorizations of $\p$.
\end{lem}

\begin{proof}
Let $f'=(\t_1',\ldots,\t_r')$ be an $r$-tuple of star transpositions that satisfies the conditions described by
Lemmas~\ref{lem:order} and~\ref{lem:between}. For brevity we shall refer to these conditions as C1 and C2,
respectively.  Suppose the cycles of $\p$ are $\s_1,\ldots,\s_m$, with $\s_1$ containing symbol 1. We wish to show
$\p' = \p$, where $\p' := \t_1' \cdots \t_r'$.  (Note that the transitivity and minimality of $f'$ are then
immediately implied by C1.)

If $\p$ has only one cycle, say $\p = (1\,b_2\,\cdots\,b_n)$, then C1 implies $\p' = \t_1' \cdots \t_r' =
(1\,b_n)(1\,b_{n-1})\cdots(1\,b_2) = \s_1$.  Hence $\p = \p'$ in this case. Otherwise, by C2 there exists some
cycle $\s_j=(a_1\,\cdots\,a_k) \neq \s_1$ of $\p$ such that the factors $\t_i'$ that meet $\s_j$ occur
contiguously in $f'$.  By C1 this means that for some $s$ we have $$\t_s' \t_{s+1}' \cdots
\t_{s+k}'=(1\,a_1)(1\,a_k)\cdots(1\,a_1)=(a_1\,a_2\,\cdots\,a_k),$$ and no factors of $f'$ other than
$\t_s',\ldots,\t_{s+k}'$ meet $\s_j$.  Thus $\p'$ agrees with $\p$ on $S:=\{a_1,\ldots,a_k\}$, and $f'' =
(\t_1',\ldots,\t_{s-1}',\t_{s+k+1}',\ldots,\t_r')$ is a star factorization of $\p'':=\p' |_{[n] \wo S}$. But $f''$
satisfies C1 and C2 relative to the permutation $\p|_{[n] \wo S}$, so we can iterate this argument to see that
$\p'$ agrees with $\p$ on all of $[n]$. \qed
\end{proof}

\section{Counting Star Factorizations}

\label{sec:counting}

\newcommand{\factor}[3]{\underset{\rule{0mm}{5mm}\displaystyle#3}{(1\,#2)}}
\newcommand{\sfactor}[3]{\underset{\rule{0mm}{5mm}\displaystyle \mathbf{\underline{#3}}}{(1\,#2)}}

Let $\p \in \Sym{n}$ be a permutation with cycles $\s_1,\ldots,\s_m$, listed in increasing order of least element
(in particular, $\s_1$ contains symbol 1). Set $r:=n+m-2$, and let $f = (\t_1,\ldots,\t_{r})$ be a minimal
transitive star factorization of $\p$. Define the word $w = w_1 \cdots w_{r} \in [m]^{r}$ by setting $w_i =j$ if
$\t_i$ meets $\s_j$. Moreover, for $2 \leq j \leq m$, define $k_j$ by the condition that the rightmost factor of
$f$ meeting $\s_j$ is $(1\,k_j)$.

\begin{exmp}
\label{exmp:word} Consider $f$ and $\p$ as defined in~\eqref{eq:fexmp} and~\eqref{eq:pexmp}. Under each factor
$\t_i$, we write the unique value of $j$ such that $\t_i$ meets $\s_j$, and we distinguish the rightmost
occurrence of each symbol $j \geq 2$:
\begin{align*}
        \factor{1}{9}{5}
        \factor{2}{\,11}{5}
        \sfactor{3}{9}{5}
        \factor{4}{2}{1}
        \factor{5}{\,10}{3}
        \factor{6}{5}{3}
        \factor{7}{3}{2}
        \sfactor{8}{3}{2}
        \factor{9}{4}{3}
        \factor{10}{7}{3}
        \factor{11}{6}{4}
        \sfactor{12}{6}{4}
        \sfactor{13}{\,10}{3}
        \factor{14}{8}{1}
\end{align*}
This yields the word
\begin{equation}
\label{eq:word}
    w = 5\;5\;5\;1\;3\;3\;2\;2\;3\;3\;4\;4\;3\;1,
\end{equation}
while the transpositions in the distinguished positions give the values of $k_j$, in this case
$(k_2,k_3,k_4,k_5)=(3,10,6,9)$. \qed
\end{exmp}

Let $\Wset{\p} \subset [m]^{r}$ be the set of words such that
\begin{itemize}
    \setlength{\itemsep}{0pt}
    \item   1 appears $\ell_1-1$ times,
    \item   $j$ appears $\ell_j+1$ times for $2 \leq j \leq m$, and
    \item   there are no occurrences of the subwords $abab$ or $a1a$ for distinct $a,b \neq 1$.
\end{itemize}
If $\orb{1},\ldots,\orb{m}$ are the orbits of $\p$, listed in increasing order of least element, then the
correspondence described above is clearly one-one between tuples $(w,k_2,\ldots,k_m) \in \Wset{\p} \times \orb{2}
\times \cdots \times \orb{m}$ and star factorizations of $\p$ satisfying the conditions of Lemmas~\ref{lem:order}
and~\ref{lem:between}. Lemma~\ref{lem:characterize} then establishes that this is, in fact, a bijection between
such tuples and the set $\Fset{\p}$ of all minimal transitive star factorizations of $\p$.  Thus we have
\begin{equation}
\label{eq:corr}
        |\Fset{\p}| = |\Wset{\p}| \cdot |\orb{2}| \cdots |\orb{m}|.
\end{equation}
We are now in a position to prove Theorem~\ref{thm:mainthm}.


\newcommand{\seq}[2]{\ensuremath{#1_0, #1_1, \ldots, #1_{#2}}}
\newcommand{\dec}[2]{#1^{(#2)}}

\noindent \textbf{Proof of Theorem~\ref{thm:mainthm}:} Assume the notation above, and for convenience let $\ell_j
= |\orb{j}|$, for $j=1,\ldots,m$. Consider the set of sequences $(\seq{d}{r})$ whose entries $d_i$ are either 1 or
$-\dec{\ell_{j}}{j}$ for some $j\geq 2$, where the exponent $(j)$ is considered to be a decoration. Of these, let
$\Dset{\pi}$ be the subset satisfying the following properties:
\begin{itemize}
    \item $-\dec{\ell_j}{j}$ appears exactly once, for $2 \leq j \leq m$, and
    \item all partial sums (ignoring decorations) are positive.
\end{itemize}
Thus $\Dset{\pi}$ describes a type of decorated Dyck sequence.

Define $\map{\Psi_\pi}{\Wset{\pi}}{\Dset{\pi}}$ by the following rule.  The image $(\seq{d}{r})$ of the word
$w_1\cdots w_r \in \Wset{\pi}$ is given by
\begin{enumerate}
    \item   $d_0 = 1$,
    \item   if $w_i=1$ then $d_i = 1$,
    \item   if $w_i=j$, for $j \geq 2$, then
                \begin{enumerate}
                \item   $d_i = -\dec{\ell_j}{j}$ when $w_{i}$
                        is the rightmost occurrence of $j$ in $w$
                \item   $d_i = 1$ otherwise.
                \end{enumerate}
\end{enumerate}
For example, with $\p$ defined by~\eqref{eq:pexmp} and $w$ given by~\eqref{eq:word}, we have
\begin{equation}
\label{eq:dexmp}
    \Psi_{\pi}(w)=(1, 1, 1, -\dec{2}{5}, 1, 1, 1, 1, -\dec{1}{2}, 1, 1, 1, -\dec{1}{4}, -\dec{4}{3}, 1).
\end{equation}

Clearly $\Psi_\pi$ is well defined.  It is also easily seen to be bijective.  Indeed, suppose ${\mathbf d} =
(\seq{d}{r}) \in \Dset{\pi}$, and let $d_i=-\dec{\ell_j}{j}$ be the first negative entry of $\mathbf{d}$.
Note that the value of $j$ is known via the decoration.  Remove $d_i$ and the previous $\ell_j$ 1s from
$\mathbf{d}$ to obtain a new sequence $\mathbf{d}'$. Inductively, $\mathbf{d}'=\Psi_{\pi'}(w')$  for a unique
word $w' \in \Wset{\pi^\prime}$, where $\pi'$ is the permutation obtained from $\pi$ by removing cycle
$\s_j$.  Adding $\ell_j + 1$ copies of $j$ after $w'_{i - \ell_j -1}$ gives the unique desired preimage
$w=\Psi_{\pi}^{-1}(\{\mathbf{d}\})$.

Thus we have $|\Wset{\pi}|=|\Dset{\pi}|$. So we now turn to enumerating $\Dset{\pi}$.  Our main tool is the cycle
lemma of Dvoretzky and Motzkin~\cite{dmot}, one version of which states that any sequence  with integral entries
$\leq 1$ and total sum $s \geq 0$ has exactly $s$ cyclic rotations with all partial sums positive.

Any sequence in $\Dset{\pi}$ has terms $-\dec{\ell_2}{2},\ldots,-\dec{\ell_m}{m}$ along with $r+1-(m-1)=n$ entries
equal to 1. Note that there are  $(n+m-1)!/n!$ sequences with exactly these terms. The sequences $(\seq{d}{r}) \in
\Dset{\pi}$ we wish to count are characterized by having total sum (ignoring decorations)
$$
    \sum_{i=0}^r d_i = n\cdot 1 - (\ell_2 + \cdots + \ell_m) = n-(n-\ell_1)=\ell_1
$$
with all partial sums positive.  Since a sequence of length $n+m-1$ admits $n+m-1$ cyclic rotations, the cycle
lemma implies that
$$
        |\Dset{\pi}| = \frac{(n+m-1)!}{n!} \cdot \frac{\ell_1}{n+m-1}.
$$
Theorem~\ref{thm:mainthm} now follows from identity~\eqref{eq:corr}, since $|\Wset{\pi}| = |\Dset{\pi}|$ and
$|\orb{j}|=\ell_j$. \qed

\newcommand{\Words}{\mathcal{W}}
\newcommand{\Restricted}{\mathcal{U}}

\section{A Graphical Approach}
\label{sec:graphical}

Transitive factorizations in the symmetric group are well known to be in correspondence with certain classes of
labelled maps, and our characterization of star factorizations (Lemmas~\ref{lem:order}, \ref{lem:between}, and
\ref{lem:characterize}) can be derived elegantly through this connection.  We now briefly describe how this is
done, using a version of the factorization-map correspondence introduced in~\cite{irving}. Indeed, it was by this
method that Theorem~\ref{thm:mainthm} was originally discovered. We elected to frame our proof in Pak's techniques
to demonstrate how they generalize and to keep this paper self contained. We note that an alternative formulation
of the factorization-map correspondence, developed with great effect in~\cite{schaeffer}, can be applied here with
equal ease.

Let $f=(\t_1,\ldots,\t_r)$ be a transitive factorization of $\p \in \Sym{n}$, where the factors $\t_i$ are
arbitrary transpositions.  Then $f$ naturally induces a graph $G_f$ on $n$ labelled vertices and $r$ labelled
edges, as follows: the vertex set of $G_f$ is $[n]$, and there is an edge with label $i$ between vertices $a$ and
$b$ whenever $\t_i=(a\,b)$. The transitivity of $f$ ensures $G_f$ is connected, so $G_f$ admits a 2-cell embedding
in an orientable surface of minimal genus. A unique such map $M_f$ is determined by insisting that the edge labels
encountered on anticlockwise traversals of small circles around the vertices are cyclically increasing.

\begin{figure}
\begin{center}
\includegraphics[width=.35\textwidth]{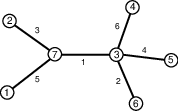}
\end{center}
\caption{The planar map corresponding to $(1\,2\,3\,4\,5\,6\,7)=(3\,7)(3\,6)(2\,7)(3\,5)(1\,7)(3\,4)$.}
\label{fig:tree}
\end{figure}

\begin{exmp}
The factorization $(1\,2\,3\,4\,5\,6\,7)=(2\,5)(3\,6)(2\,7)(3\,5)(1\,7)(3\,4)$ is minimal transitive.  Its
corresponding planar map is shown in Figure~\ref{fig:tree}.  \qed
\end{exmp}

As described in~\cite{irving}, faces of $M_f$ correspond with the cycles of $\p$.  In particular, let $F$ be a
face of $M_f$, and let $(e_0,\ldots,e_m)$ be the cyclic list of edge labels encountered along a counterclockwise
traversal of the boundary of $F$. If $i_1,\ldots,i_k$ index the \emph{ascents} of this list (that is, $e_i \leq
e_{i+1}$ if and only if $i \in \{i_1,\ldots,i_k\}$), then $\pi$ contains the cycle $(a_1\,a_2\,\cdots\,a_k)$,
where $a_j$ is the label of the vertex at the corner of $M_f$ formed by edges $e_{i_j}$ and $e_{i_j+1}$.

\begin{figure}
\begin{center}
\includegraphics[width=.95\textwidth]{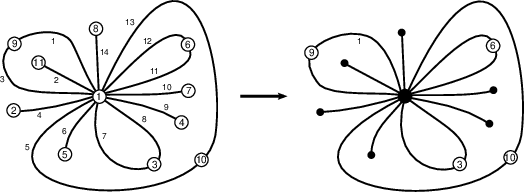}
\end{center}
\caption{The planar map corresponding to a star factorization, and its reduced form.} \label{fig:map}
\end{figure}

With this correspondence, the Euler-Poincar\'e formula implies \emph{$M_f$ is planar precisely when $f$ is minimal
transitive.} Indeed, the maps corresponding to minimal transitive \emph{star} factorizations are particularly
simple. This is illustrated in Figure~\ref{fig:map}, where the planar map associated with our primary example
factorization~\eqref{eq:fexmp} is drawn.

Since such a map must be planar with edge labels increasing clockwise around the central vertex 1, no edge
$\{1,a\}$ can appear more than twice. When two copies of $\{1,a\}$ are present they enclose a face of the map. It
is this face that is associated with the cycle of the target permutation containing symbol $a$, and a vertex $b$
of degree one lies within it precisely when $b$ belongs to this same cycle.  Translated from the language of maps
to that of factorizations, these observations are equivalent to Lemmas~\ref{lem:order} and~\ref{lem:between}.

The canonical labelling of edges around the central vertex makes all but label 1 superfluous. Moreover, the labels
of all vertices of degree 1 may be deduced from the target permutation and the labels of the other vertices.  Thus
all maps corresponding to minimal transitive star factorizations may be reduced in the manner demonstrated on the
right of Figure~\ref{fig:map}.

\begin{figure}
\begin{center}
\includegraphics[width=\textwidth]{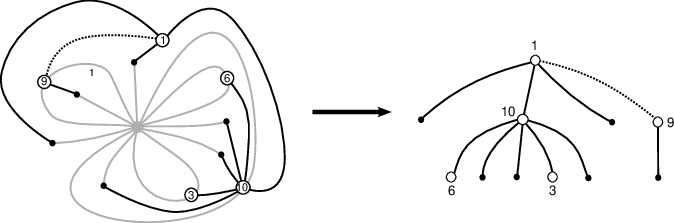}
\end{center}
\caption{From maps to  trees.} \label{fig:maptotree}
\end{figure}

From this reduced form, create a rooted plane tree as follows.  Begin by placing a root vertex with label 1 in the
outer face. Every labelled vertex is now naturally associated with one face of the map.  Then draw an edge between
each labelled vertex and all (non-central) vertices lying within with its associated face.  One of these edges
will join vertex 1 to the endpoint of the map edge with label 1. This is to be considered the root edge of the
tree.  See Figure~\ref{fig:maptotree} for an example.

This transition from factorization to map to tree is reversible. A minimal transitive factorization of a
permutation $\p \in \Sym{n}$ with orbits $\orb{1},\ldots,\orb{m}$ (listed, as usual, in increasing order of least
element) corresponds with a tree on $m$  labelled white vertices and $n-m$ black vertices in which
\begin{enumerate}
\item   the root is white with label 1,
\item   the non-root white vertices are labelled $\{a_2,\ldots,a_m\}$, where $a_j \in \orb{j}$,
\item   the white vertex with label $a_j$ has $|\orb{j}|-1$ black children, for $j=1,\ldots,m$.
\end{enumerate}
Such trees can be encoded using Dyck-type sequences, as follows:  traverse the boundary, beginning at the root and
proceeding clockwise along the root edge, writing 1 whenever a vertex is encountered for the first time, and
$-\dec{i}{j}$ when a white vertex with label $j \geq 2$ and $i-1$ black children is encountered for the last time.
For instance, the tree in Figure~\ref{fig:maptotree} yields the following sequence (compare
with~\eqref{eq:dexmp}):
$$(1,1,1,\dec{-2}{9},1,1,1,1,\dec{-1}{3},1,1,1,\dec{-1}{6},\dec{-4}{10},1).$$ These sequences are counted as
in Section~\ref{sec:counting} to yield Theorem~\ref{thm:mainthm}.

\section{Further Questions}
\label{sec:final}

Notice that Theorem \ref{thm:mainthm} asserts that the number of minimal transitive star factorizations of a
permutation $\pi$ depends only on the conjugacy class of $\pi$ (that is, the length of its cycles).  This is not
obvious from the formulation of the problem, since one would certainly expect that the length of the cycle of
$\pi$ containing symbol 1 would play a special role.

Moreover, while this article was being refereed, Goulden and Jackson~\cite{gj-jucys} extended
Theorem~\ref{thm:mainthm} to compute the number of transitive star factorizations of any permutation into an
\emph{arbitrary} number of factors (that is, minimality is not assumed).  Interestingly, they witness the same
symmetry in their results: the number of transitive star factorizations of $\pi$ of length $r$ is dependent only
on the conjugacy class of $\pi$.

Finding a simple combinatorial explanation for this curious symmetry remains an interesting open problem.  Further
open questions regarding star factorizations and their role in the general interplay between factorizations and
geometry are discussed in~\cite{gj-jucys}.

\section*{Acknowledgements}

Both authors would like to thank Ian Goulden for some useful discussions, and an anonymous referee for a number of
helpful comments.  AR is supported by a Natural Sciences and Engineering Research Council of Canada Postdoctoral
Fellowship.


\end{document}